\documentclass[12pt]{article}

\usepackage{graphicx,tikz,color,url}
\usepackage{amsmath,amssymb,latexsym}
\usetikzlibrary{arrows}

\newtheorem{theorem}{Theorem}[section]

\newtheorem{lemma}[theorem]{Lemma}
\newtheorem{proposition}[theorem]{Proposition}

\newtheorem{corollary}[theorem]{Corollary}

\newcommand{\proof}{\noindent{\bf Proof.\ }}
\newcommand{\qed}{\hfill $\square$ \bigskip}

\newcommand{\gsmb}{\gamma_{\rm SMB}}
\newcommand{\gmb}{\gamma_{\rm MB}}

\begin{document}

\title{Maker-Breaker domination game played on corona products of graphs}

\author{Athira Divakaran$^{a,}$\thanks{Email: \texttt{athiradivakaran91@gmail.com}} 
\and 
Tijo James$^{b,}$\thanks{Email: \texttt{tijojames@gmail.com}} 
\and 
Sandi Klav\v zar $^{c,d,e,}$\thanks{Corresponding author. Email: \texttt{sandi.klavzar@fmf.uni-lj.si}}
\and 
Latha S Nair$^{a,}$\thanks{Email: \texttt{lathavichattu@gmail.com}} 
}
\maketitle

\begin{center}
$^a$ Department of Mathematics, Mar Athanasius College, Kothamangalam, India\\
	\medskip
$^b$ Department of Mathematics, Pavanatma College, Murickassery, India\\
\medskip

	$^c$ Faculty of Mathematics and Physics, University of Ljubljana, Slovenia\\
	\medskip

	$^d$ Institute of Mathematics, Physics and Mechanics, Ljubljana, Slovenia\\
	\medskip
	
	$^e$ Faculty of Natural Sciences and Mathematics, University of Maribor, Slovenia\\
	\medskip
 
\end{center}

\begin{abstract}
In the Maker-Breaker domination game, Dominator and Staller play on a graph $G$ by taking turns in which each player selects a not yet played vertex of $G$. Dominator's goal is to select all the vertices in a dominating set, while Staller aims to prevent this from happening. In this paper, the game is investigated on corona products of graphs. Its outcome is determined as a function of the outcome of the game on the second factor. Staller-Maker-Breaker domination numbers are determined for arbitrary corona products, while Maker-Breaker domination numbers of corona products are bounded from both sides. All the bounds presented are demonstrated to be sharp. Corona products as well as general graphs with small (Staller-)Maker-Breaker domination numbers are described. 
\end{abstract}


\noindent
{\bf AMS Subj.\ Class.\ (2010)}: 05C57, 05C69
  
\section{Introduction}

The Maker-Breaker game was introduced by Erd\H{o}s and Selfridge in~\cite{erdos-1973}. The game is played on an arbitrary hypergraph by two players named Maker and Breaker. The players alternately select an unplayed vertex during the game. Maker's aim is to occupy all the vertices of some hyperedge, the goal of Breaker is to prevent him from doing it. The game has been extensively researched, both in general and in specific cases, cf.\ the book~\cite{hefetz-2014}, the recent paper~\cite{rahman-2023}, and references therein.  

In this paper we are interested in the domination version of the Maker-Breaker game which was introduced in 2020 by Duch\^{e}ne, Gledel, ~Parreau, and Renault~\cite{duchene-2020}. The {\em Maker-Breaker domination game} ({\em MBD game} for short) is a game played on a graph $G=(V(G), E(G))$  by two players named Dominator and Staller. These names were chosen so that the players are named in line with the previously intensively researched domination game~\cite{bresar-2010, book-2021}. Just as in the general case, the two players alternately select unplayed vertices of $G$. The aim of Dominator is to select all the vertices of some dominating set of $G$, while Staller aims to select at least one vertex from every dominating set of $G$. There are two variants of this game depending on which player has the first move. A {\em D-game} is the MBD game in which Dominator has the first move and an {\em S-game} is the MBD game in which Staller has the first move. 

The following graph invariants are naturally associated to the MBD game~\cite{bujtas-2024, gledel-2019}. The {\em Maker-Breaker domination number}, $\gmb(G)$, is the minimum number of moves of Dominator to win the D-game on $G$ when both players play optimally. That is, $\gmb(G)$ is the minimum number of moves of Dominator such that he wins in this number of moves no matter how Staller is playing. If Dominator has no winning strategy in the D-game, then set $\gmb(G)= \infty$. The {\em Staller-Maker-Breaker domination number},  $\gsmb(G)$, is the minimum number of moves of Staller to win the D-game on $G$ when both players play optimally, where $\gsmb(G)= \infty$ if Staller has no winning strategy.  In a similar manner, $\gmb'(G)$ and $\gsmb'(G)$ are the two parameters associated with  the S-game. Briefly, we will call $\gmb(G)$, $\gmb'(G)$, $\gsmb(G)$, and $\gsmb'(G)$ the {\em MBD numbers of} $G$.

In the seminal paper~\cite{duchene-2020} it was proved, among other things, that deciding the winner of the MBD game can be solved efficiently on trees, but it is PSPACE-complete even for bipartite graphs and split graphs. In~\cite{duchene-2023+} the authors give a complex linear algorithm for the Maker-Maker version of the game. The paper~\cite{bujtas-2024} focuses on $\gsmb(G)$ and $\gsmb'(G)$ and among other results establishes appealing exact formulas for $\gsmb'(G)$ where $G$ is a path. In~\cite{bujtas-2023}, for every positive integer $k$, trees $T$ with $\gsmb(T) = k$ are characterized and exact formulas for $\gsmb(G)$ and $\gsmb'(G)$ derived for caterpillars. In the main result of~\cite{forcan-2023+}, $\gmb$ and $\gmb'$ are determined for Cartesian products of $K_2$ by a path. In~\cite{dokyeesun-2024+}, the MBD game is further studied on Cartesian products of paths, stars, and complete bipartite graphs. The total version of the MBD game was introduced in~\cite{gledel-2020} and further investigated in~\cite{forcan-2022}.

The corona product of graphs is a graph operation that has already been studied from many perspectives, \cite{alagammai-2020, dong-2021, klavzar-2022, nie-2023, nirajan-2023, verma-2024} is a selection of recent relevant contributions. In this paper, we investigate the MBD game played on these products and proceed as follows. In the rest of the introduction definitions and known results needed are stated. In the next section we determine the outcome of the MBD game on an arbitrary corona product of graphs. In Section~\ref{sec:MBD-number} we determine $\gsmb$ and $\gsmb'$, and bound $\gmb$ and $\gmb'$ of arbitrary corona products. We demonstrate that all the bounds are sharp. In particular, $\gmb$ and $\gmb'$ are determined for corona products of arbitrary graphs by trees with perfect matchings. In Section~\ref{sec:small MB and SMB values} we describe corona products with small MBD numbers and along the way  describe general graphs with small MBD numbers.  

Let $G=(V(G), E(G))$ be a  graph. The order of  $G$ will be denoted by $n(G)$. For a vertex $v\in V(G)$,  its open neighbourhood is denoted by $N(v)$ and its closed neighbourhood by $N[v]$.  The degree of $v$  is $\deg(v)=|N(v)|$.  The minimum and the maximum degree of $G$ are denoted  by $\delta(G)$ and $\Delta(G)$. An isolated vertex is a vertex of degree $0$, a leaf is a vertex of degree $1$.  A {\em support vertex} is a vertex adjacent to a leaf, a  {\em strong support vertex} is a vertex that is adjacent to at least two leaves. A {\em dominating set} of $G$ is a set $D \subseteq V(G)$ such that each vertex from $V(G)\setminus D$ has a neighbour in $D$. The {\em domination number} $\gamma(G)$ of $G$  is the minimum cardinality of a dominating set of $G$. A dominating set of $G$ of cardinality $\gamma(G)$ is refered to as a {\em $\gamma$-set} of $G$. A vertex of degree $n(G) -1$ is a {\em dominating vertex}.  

The \textit{corona} $G\odot H$  of graphs $G$ and $H$ is obtained by taking one copy of $G$ and $n(G)$ copies of $H$ by joining the $i^{{\rm th}}$ vertex of $G$ to each vertex in the $i^{{\rm th}}$ copy of $H$. Clearly, $n(G\odot H)=n(G)(n(H)+1)$ and $\delta(G\odot H)=\delta(H)+1$. It is also straightforward to see that $\gamma(G\odot H) = n(G)$.  

We next collect known results about the MBD game to be used later on. 

\begin{lemma} [No-Skip Lemma for Dominator~\cite{gledel-2019}] 
\label{22} 
In an optimal strategy of Dominator to achieve $\gmb(G)$ or $\gmb'(G)$ it is never an advantage for him to skip a move. Moreover, if Staller skips a move it can never disadvantage Dominator.  
\end{lemma}

\begin{lemma} [No-Skip Lemma for Staller~\cite{bujtas-2024}] 
\label{22 A} 
In an optimal strategy of Staller to achieve $\gsmb(G)$ or $\gsmb'(G)$ it is never an advantage for her to skip a move. Moreover, if Dominator skips a move it can never disadvantage Staller.  
\end{lemma}

For $X\subseteq V(G)$, let $G|X$ denote the graph $G$ in which vertices from $X$ are considered as being already dominated.

\begin{theorem} [Continuation Principle~\cite{gledel-2019}]
\label{thm: Continuation Principle}
Let $G$ be a graph with $A, B \subseteq V(G)$. If $B\subseteq A$ then $\gmb(G|A) \leq \gmb(G|B)$ and  $\gmb'(G|A) \leq \gmb'(G|B)$.
\end{theorem}

\begin{theorem} {\rm \cite{bujtas-2024}}
\label{22F}
Let G be a graph.
\begin{enumerate}
\item If $\gsmb'(G) < \infty$, then $\gsmb'(G) \leq \left \lceil\frac{n(G)}{2}\right \rceil$.
\item If $\gsmb(G) <\infty$, then $\gsmb(G) \leq \left \lfloor\frac{n(G)}{2}\right \rfloor$.
\end{enumerate}
Moreover, both bounds are sharp.
 \end{theorem}
 
\section{Outcome of the game on corona products}
\label{sec:outcome}

In this section we determine the winner of the MBD game played on corona products. Before that we argue that in the rest we can restrict our attention to corona products in which the first factor is of order at least two. 

Let $H$ be an arbitrary graph. Then it is straightforward to verify that 
\begin{align*}
\gmb(K_1\odot H ) & = 1\,, \\
\gmb'(K_1\odot H ) & = \gmb(H)\,, \\
\gsmb(K_1\odot H ) & = \infty\,, \\
\gsmb'(K_1\odot H ) & = \gsmb(H)\,.
\end{align*}
Therefore, we may indeed restrict in the rest of the paper to corona products $G\odot H$, where $G$ is a graph with $n(G)\ge 2$, where $n(G)$ denotes the order of $G$. In what follows, we will also always assume that $V(G)=\{v_1,\ldots, v_{n(G)}\}\subseteq V(G\odot H)$, and denote by $H_i$ the copy of $H$ in $G\odot H$ assigned to $v_i$, $i\in [n(G)]$. 

It was proved in~\cite{duchene-2020} that for the {\em outcome} $o(G)$ of the MBD game played on $G$ we have $o(G)\in \{\mathcal{D}, \mathcal{S}, \mathcal{N}\}$, with the following meaning: 
     \begin{itemize}
         \item $o(G)=\mathcal{D}$: Dominator has a winning strategy no matter who starts the game; 
         \item $o(G)=\mathcal{S}$: Staller has a winning strategy no matter who starts the game;
         \item  $o(G)=\mathcal{N}$: the first player has a winning strategy.
    \end{itemize}

We now prove that the winner in the MBD game played on a corona product depends only on the winner of the game played on the second factor. 

\begin{theorem}\label{thm:outcome}
If $G$ is a graph with $n(G)\geq 2$ and $H$ is a graph, then 
$$o(G\odot H) = 
\begin{cases}
\mathcal{D};&  o(H)\in \{\mathcal{D},\mathcal{N}\},\\
\mathcal{S}; &  o(H)=\mathcal{S}. 
\end{cases}$$ 
\end{theorem}

\proof
Suppose first that $o(H)\in \{\mathcal{D},\mathcal{N}\}$. Then in both cases we have $\gmb(H)<\infty$. We describe a strategy for Dominator by which he can win the MBD game played on $ G\odot H$, no matter who is the first player. First, if Staller selects a vertex $v_i$ as her first move, then Dominator replies by a vertex in $H_i$ which is his optimal first move when the game is played on $H$. During later stages of the game, if Staller at some point selects another vertex of $H_i$, Dominator replies by an optimal move (w.r.t.\ the game played in $H$) in $H_i$. Playing like that, at the end of the game Dominator will dominate the vertices of $H_i$. Second, if Staller first selects some vertex of $H_i$, then Dominator replies by playing $v_i$ and then clearly at the end of the game he will dominate the vertices of $H_i$. Following the described strategy of Dominator we conclude that $o(G\odot H) = \mathcal{D}$ when $o(H)\in \{\mathcal{D},\mathcal{N}\}$. 

Suppose second that  $o(H)=\mathcal{S}$. As $n(G)\ge 2$, Staller can select some vertex $v_i$ as her first move, no matter who starts the game, such that no vertex of $H_i$ has yet been played. In the continuation of the game, Staller can isolate a vertex in $H_i$ because $o(H)=\mathcal{S}$, that is, she follows her optimal strategy on $H_i$. We conclude that $o(H)=\mathcal{S}$. 
\qed

\section{MBD numbers of corona products}
\label{sec:MBD-number}

In this section we consider the MBD numbers of corona products $G\odot H$. In view of Theorem~\ref{thm:outcome} we consider three cases depending on $o(H)$. We start with the simplest case $o(H)=\mathcal{S}$ which enables an exact solution as follows. 

\begin{theorem}\label{thm: formula 3}
If $o(H)=\mathcal{S}$ and $G$ is a graph with $n(G) \geq  2$, then  $$\gsmb(G\odot H)=\gsmb'(G\odot H)=1+\gsmb(H).$$
\end{theorem}

\proof
Using Theorem~\ref{thm: Continuation Principle}, some vertex of $G$ is an optimal first move for Dominator, say $v_i$. Then all the vertices in its closed neighbourhood in $G\odot H$ are dominated. Staller responds by selecting a vertex $v_j$, where $j \neq i$. Staller's plan afterwards is to play only on $H_j$. If Dominator plays his second move on $H_j$, then Staller can finish the game by selecting $\gsmb(H)$ vertices of $H_j$ and if Staller is able to play first in $H_j$, this can be done in $\gsmb'(H)$ moves by Staller. As $\gsmb(H) \ge \gsmb'(H)$, see~\cite[Proposition 2.7]{bujtas-2024}, Staller thus has a strategy to end the game in at most $1 + \gsmb(H)$ moves. On the other hand, Dominator can force this number of moves by playing (optimally in the S-game on $H$) his second move in $H_j$, hence $\gsmb(G\odot H) \le 1+\gsmb(H)$. 

The above argument is actually independent of who moves first since $n(G)\ge 2$ and Staller can select a vertex of $G$ as her first move. Therefore we also have $ \gsmb'(G\odot H)= 1+ \gsmb(H)$.
\qed

The next case we consider is $o(H)=\mathcal{D}$, for which the following holds. 

\begin{theorem}
\label{thm:bounds 1}
If $o(H)=\mathcal{D}$ and $G$ is a graph with $n(G)\geq 2$, then   
$$\left \lceil \frac{n(G)}{2} \right \rceil +   \left \lfloor \frac{n(G)}{2} \right \rfloor \gmb(H) \le \gmb(G\odot H) \le \left \lceil \frac{n(G)}{2} \right \rceil +   \left \lfloor \frac{n(G)}{2} \right \rfloor \gmb'(H)$$  
and  
$$\left \lfloor \frac{n(G)}{2} \right \rfloor +  \left \lceil \frac{n(G)}{2} \right \rceil\gmb(H) \le \gmb'(G\odot H) \le \left \lfloor \frac{n(G)}{2} \right \rfloor +  \left \lceil \frac{n(G)}{2} \right \rceil\gmb'(H)\,.$$
\end{theorem}

\proof
Consider the following strategy of Dominator in the D-game played on $G\odot H$. His basic idea is to select as many vertices of $G$ as possible. In particular, if at some point of the game Staller plays a vertex of $H_i$ and the vertex $v_i$ has not yet been played, then Dominator selects $v_i$ as his next move. On the other hand, if at some point of the game Staller plays a vertex $w$ of $H_i$ and the vertex $v_i$ has been already selected by Staller, then Dominator replies by a vertex in $H_i$ which is an optimal reply for him in the game played on $H\cong H_i$ after the first move $w$ of Staller. In the continuation of the game, whenever Staller plays in some $H_i$, Dominator replies optimally there. It is also possible that Dominator is the first to play in some $H_i$ in which no vertex has yet been played (but $v_i$ was played by Staller), in which case Dominator selects a vertex which is optimal in the D-game played on $H$. 

Following the above strategy, Dominator will select at least $\left \lceil \frac{n(G)}{2} \right \rceil$ vertices of $G$. If Staller was the first to select a vertex in $H_i$, then Dominator will select at most at most $\gmb'(H)$ vertices of $H_i$. Similarly, If Dominator was the first to select a vertex in $H_i$, then he will select at most at most $\gmb(H)$ vertices of $H_i$. Since $\gmb(H) \le \gmb'(H)$ (see~\cite[Theorem 3.1]{gledel-2019}), Dominator will play at most $\gmb'(H)$ vertices of each $H_i$, where $v_i$ was selected by Staller. As there are at most $\left \lfloor \frac{n(G)}{2} \right \rfloor$ such vertices, the upper bound for $\gmb(G\odot H)$ follows. 

To prove the lower bound for $\gmb(G\odot H)$, consider the following strategy of Staller. During the first phase of the game, Staller selects as many vertices of $G$ as possible, that is, $\left \lfloor \frac{n(G)}{2} \right \rfloor$ vertices. After that, whenever Dominator plays a vertex of $H_i$, Staller replies there by her optimal strategy in the game played on $H$. In this way Staller forces to select at least $\gmb(H)$ vertices in $H_i$ if Dominator was the first to play in $H_i$ and at least $\gmb'(H)$ vertices if Staller was the first to play in $H_i$. Since $\gmb'(H) \ge \gmb(H)$, in each of the $\left \lfloor \frac{n(G)}{2} \right \rfloor$ graphs $H_i$ for which $v_i$ was played by Staller, at least $\gmb(H)$ vertices will be played by Dominator. This proves the lower bound. 

The arguments for the bounds for $\gmb'(G\odot H)$ are parallel and hence omitted. 
\qed

The bounds of Theorem~\ref{thm:bounds 1} are sharp as can be demonstrated by the following two consequences of the theorem, where the first one directly follows from it. 

\begin{corollary}
\label{cor:formula 1}
If $o(H)= \mathcal{D}$ with $\gmb(H)=\gmb'(H)$, and $G$ is a graph with $n(G)\geq 2$, then   
$$\gmb(G\odot H)=\left \lceil \frac{n(G)}{2} \right \rceil +   \left \lfloor \frac{n(G)}{2} \right \rfloor \gmb(H)$$  
and  
$$\gmb'(G\odot H)=\left \lfloor \frac{n(G)}{2} \right \rfloor +  \left \lceil \frac{n(G)}{2} \right \rceil\gmb(H).$$
\end{corollary}

\begin{theorem}\label{33}
Let $G$ be a graph with $n(G)\ge 2$ and let $T$ be a tree with a perfect matching. Then  
$$\gmb(G\odot T)=\left \lceil \frac{n(G)}{2} \right \rceil +   \left \lfloor \frac{n(G)}{2} \right \rfloor \frac{n(T)}{2}$$ 
and  
$$\gmb'(G\odot T)=\left \lfloor \frac{n(G)}{2} \right \rfloor +  \left \lceil \frac{n(G)}{2} \right \rceil \frac{n(T)}{2}\,.$$
\end{theorem}

\proof
Since $T$ is a tree with a perfect matching, we know from~\cite{duchene-2020} that $o(T) = \mathcal{D}$ as well as that $\gmb(T) = \gmb'(T) = \frac{n(T)}{2}$. The conclusion then follows by Corollary~\ref{cor:formula 1}. 
\qed

The last case to consider is $o(H)=\mathcal{N}$, for which the following holds. 

\begin{theorem}
\label{thm: formula 2}
If $o(H)=\mathcal{N}$ and $G$ is a graph with $n(G)\geq 2$, then 
$$1 + \left\lfloor \frac{n(G)-1}{2} \right\rfloor \gamma(H) +   \left \lfloor \frac{n(G)}{2}\right \rfloor \gmb(H) \le \gmb(G\odot H) \le 1 +(n(G)-1)\gmb(H)$$ 
and
$$ \left\lfloor \frac{n(G)}{2} \right\rfloor \gamma(H) +   \left \lceil \frac{n(G)}{2} \right\rceil \gmb(H) \le \gmb'(G\odot H) \le n(G)\gmb(H)\,.$$ 
\end{theorem}

\proof
Consider the D-game on $G\odot H$. 

To prove the lower bound, consider the following strategy of Staller. After Dominator selects some vertex of $G$, say $v_i$, as his optimal first move in view of Theorem~\ref{thm: Continuation Principle}, Staller selects as her first move another vertex of $G$, say $v_j$. Because $o(H) = \mathcal{N}$, Dominator must reply by playing (an optimal move with respect to the game played on $H$) in $H_j$. After that Staller continues by selecting the vertices of $G$ until all are vertices of $G$ but $v_i$ are played by Staller. For each such vertex $v_k$ Dominator was forced to select a vertex from $H_k$. Hence until this point of the game Dominator played $n(G)$ vertices, the vertex $v_i$ and one vertex in each $H_j$, $j\ne i$. In the rest of the game, Staller will be able to play the second move in at least $\left \lfloor \frac{n(G)}{2}\right \rfloor$ copies $H_j$ of $H$. Consequently, she will be able to force at least $\gmb(H)$ moves of Dominator in each of these copies. As in the remaining $\left\lfloor \frac{n(G)-1}{2} \right \rfloor$ copies of $H$ at least $\gamma(H)$ moves must be done  Dominator, the bound follows. 

For the upper bound, we argue as follows. After Dominator selects a vertex $v_i$ as his first move, Staller responds to this move by playing a vertex $v_j$, $j\neq i$. Otherwise, if Staller selects a vertex  in  $H_j$ then Dominator can dominate all the vertices of $H_j$ by selecting the vertex  $v_j$ in a single move. Since $\gmb(H) < \infty$,  Dominator has an optimal strategy in the D-game on $H$  to finish it in at most $\gmb(H)$ moves. Dominator must play a vertex of  $H_j$ as his response after the Staller's move on $v_j$ because $\gmb'(H) =\infty $. Now Dominator's strategy is as follows. If Staller plays on $H_j$ then Dominator uses his optimal strategy on $H_j$ for the remaining moves in $H_j$. Otherwise, if Staller selects a vertex of $G$, say $v_k$, then Dominator uses his optimal strategy in $H_k$ for the remaining moves in $H_k$. In this strategy, Dominator can finish the game on $G\odot H$ in $1+(n(G)-1)\gmb(H)$ moves. 
Here, Dominator may not play optimally  but Staller plays optimally in $G \odot H$. Therefore $\gmb(G\odot H) \leq 1+(n(G)-1)\gmb(H)$.

Similar arguments also hold for the S-game on $G\odot H$. The only difference is that Staller starts the game and hence she can select all the vertices of $G$ in $G\odot H$. After that Dominator can ensure that he will play at most $\gmb(H)$ vertices in each $H_i$ by replying to each move of Staller in some $H_j$ by an optimal move in the same $H_j$. As for the lower bound, Staller can now force Dominator to play at least $\gmb(H)$ moves in $\left \lceil \frac{n(G)}{2} \right\rceil$ copies $H_j$ of $H$. 
\qed

All the bounds of Theorem~\ref{thm: formula 2} are sharp as follows from the next consequence of Theorem~\ref{thm: formula 2}. 

\begin{corollary}
\label{cor: all-bounds-sharp}
If $H$ is a graph with $o(H)=\mathcal{N}$ and $\gmb(H) = \gamma(H)$,  and $G$ is a graph with $n(G)\geq 2$, then 
$$\gmb(G\odot H) = 1 +(n(G)-1)\gmb(H)$$ 
and
$$\gmb'(G\odot H) = n(G)\gmb(H)\,.$$ 
In particular, if $\gmb(H) = 2$, then $\gmb(G\odot H) = 2n(G) - 1$ and $\gmb'(G\odot H) = 2n(G)$.
\end{corollary}

\proof
Since $\gmb(H) = \gamma(H)$ and $\left\lfloor \frac{n(G)-1}{2} \right\rfloor +   \left \lfloor \frac{n(G)}{2}\right \rfloor = n-1$, the bounds for $\gmb(G\odot H)$ in Theorem~\ref{thm: formula 2} coincide, and since $\left\lfloor \frac{n(G)}{2} \right\rfloor +  \left \lceil \frac{n(G)}{2} \right\rceil = n$ the same conclusion holds for the bounds for $\gmb'(G\odot H)$. 

If $\gmb(H) = 2$, then $G$ clearly does not contain a dominating vertex, hence $\gamma(H) \ge 2$ and therefore $\gamma(H) = 2$. The last assertion of the corollary now follows from the first one.  
\qed

Let $H_m$, $m\ge 2$, be the graph obtained from $m$ disjoint copies of $C_4$ by selecting a vertex in each of the copies and identifying these vertices. (The identified vertex is of degree $2m$ in $H_m$.) Then it is not hard to check that $o(H_m) = \mathcal{N}$ and that $\gmb(H_m) = \gamma(H_m) = m+1$. Hence the graphs $H_m$ fulfill the assumptions of Corollary~\ref{cor: all-bounds-sharp}.  

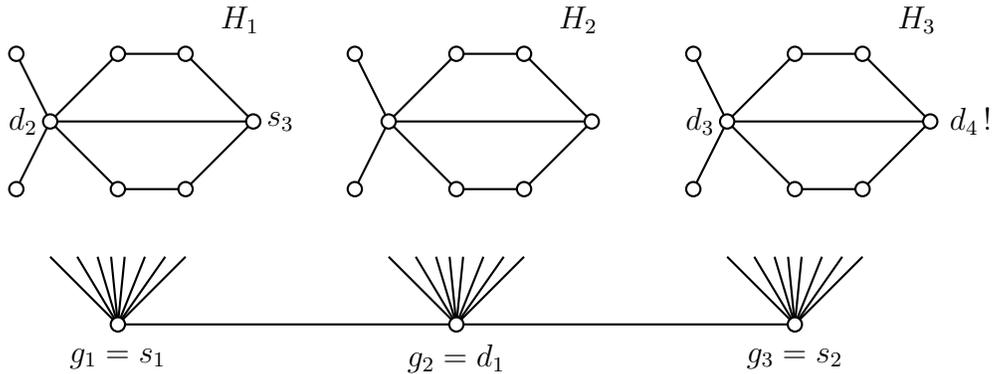
\begin{figure}[ht!]
\begin{center}
\begin{tikzpicture}[scale=0.9,style=thick,x=1cm,y=1cm]
\def\vr{3pt}

\begin{scope}[xshift=0cm, yshift=0cm] 
\coordinate(g) at (2,0);
\coordinate(h1) at (0.5,2);
\coordinate(h2) at (0.5,4);
\coordinate(h3) at (1,3);
\coordinate(h4) at (2,2);
\coordinate(h5) at (2,4);
\coordinate(h6) at (3,2);
\coordinate(h7) at (3,4);
\coordinate(h8) at (4,3);
\draw (h1) -- (h3) -- (h2); 
\draw (h3) -- (h4) -- (h6) -- (h8) -- (h7) -- (h5) -- (h3) -- (h8); 
\draw (g) -- (1,1);
\draw (g) -- (1.4,1);
\draw (g) -- (1.7,1);
\draw (g) -- (1.9,1);
\draw (g) -- (2.1,1);
\draw (g) -- (2.4,1);
\draw (g) -- (2.7,1);
\draw (g) -- (3,1);
\draw (g) -- (12,0);
\draw(g)[fill=white] circle(\vr);
\foreach \i in {1,...,8}
{ \draw(h\i)[fill=white] circle(\vr);
}
\node at (2,-0.5) {$g_1 = s_1$};
\node at (0.6,3) {$d_2$};
\node at (4.4,3) {$s_3$};
\node at (3.8,4.5) {$H_1$};
\end{scope}

\begin{scope}[xshift=5cm, yshift=0cm] 
\coordinate(g) at (2,0);
\coordinate(h1) at (0.5,2);
\coordinate(h2) at (0.5,4);
\coordinate(h3) at (1,3);
\coordinate(h4) at (2,2);
\coordinate(h5) at (2,4);
\coordinate(h6) at (3,2);
\coordinate(h7) at (3,4);
\coordinate(h8) at (4,3);
\draw (h1) -- (h3) -- (h2); 
\draw (h3) -- (h4) -- (h6) -- (h8) -- (h7) -- (h5) -- (h3) -- (h8); 
\draw (g) -- (1,1);
\draw (g) -- (1.4,1);
\draw (g) -- (1.7,1);
\draw (g) -- (1.9,1);
\draw (g) -- (2.1,1);
\draw (g) -- (2.4,1);
\draw (g) -- (2.7,1);
\draw (g) -- (3,1);
\draw(g)[fill=white] circle(\vr);
\foreach \i in {1,...,8}
{ \draw(h\i)[fill=white] circle(\vr);
}
\node at (2,-0.5) {$g_2 = d_1$};
\node at (3.8,4.5) {$H_2$};
\end{scope}

\begin{scope}[xshift=10cm, yshift=0cm] 
\coordinate(g) at (2,0);
\coordinate(h1) at (0.5,2);
\coordinate(h2) at (0.5,4);
\coordinate(h3) at (1,3);
\coordinate(h4) at (2,2);
\coordinate(h5) at (2,4);
\coordinate(h6) at (3,2);
\coordinate(h7) at (3,4);
\coordinate(h8) at (4,3);
\draw (h1) -- (h3) -- (h2); 
\draw (h3) -- (h4) -- (h6) -- (h8) -- (h7) -- (h5) -- (h3) -- (h8); 
\draw (g) -- (1,1);
\draw (g) -- (1.4,1);
\draw (g) -- (1.7,1);
\draw (g) -- (1.9,1);
\draw (g) -- (2.1,1);
\draw (g) -- (2.4,1);
\draw (g) -- (2.7,1);
\draw (g) -- (3,1);
\draw(g)[fill=white] circle(\vr);
\foreach \i in {1,...,8}
{ \draw(h\i)[fill=white] circle(\vr);
}
\node at (2,-0.5) {$g_3 = s_2$};
\node at (0.6,3) {$d_3$};
\node at (4.6,3) {$d_4$\,!};
\node at (3.8,4.5) {$H_3$};
\end{scope}

\end{tikzpicture}
\caption{An example in which it is useful for Dominator to skip a move in some copy of $H$}
	\label{fig:ilustrative-example}
\end{center}
\end{figure}

On the other hand, the bounds of Theorem~\ref{thm: formula 2} are not sharp in general. To see this, consider the example from Fig.~\ref{fig:ilustrative-example} in which $G = P_3$ and $H$ is clear from the figure. To make the figure clearer, the edges between the vertices of $G$ and $H$ are only indicated. 

Note first that $o(H) = \mathcal{N}$. in particular, the support vertex of $H$ is the optimal first move of Staller in the S-game. Moreover, the same vertex is the optimal first move of Dominator in the D-game and we have $\gmb(H) = 3$. Consider now the D-game played on $P_3\odot H$, see Fig.~\ref{fig:ilustrative-example} again. In view of Theorem~\ref{thm: Continuation Principle}, Dominator begins the game by selecting $d_1 = g_2$. After Staller replies by $s_1 = g_1$, Dominator is forced to play $d_2$ in $H_1$, and after the move $s_2 = g_3$ Dominator must play in $H_3$. Assume without loss of generality that afterwards Staller plays a vertex of $H_1$. Now it is useful for Dominator not to follow Staller in $H_1$ but instead playing the vertex $d_4$ of $H_3$ as shown in the figure. In this way, Dominator will finish the game in six moves, which is also optimal for him, that is, $\gmb(P_3 \odot H) = 6$ which coincides with the lower bound of the theorem. On the other hand, the upper bound of Theorem~\ref{thm: formula 2} yields $\gmb(P_3 \odot H) \le 7$.  

We conclude with another consequence of Theorem~\ref{thm: formula 2} which demonstrates the sharpness of bounds. 

\begin{corollary}
If $o(H)= \mathcal{D}$, then $\gmb(K_2\odot H) = 1 + \gmb(H)$.   
\end{corollary}

\proof
In this case the bounds of Theorem~\ref{thm: formula 2} coincide. 
\qed

\section{Corona products with small MBD numbers}
\label{sec:small MB and SMB values}

In this section we characterize corona products whose MBD numbers are equal to $2$ or $3$. For this purpose we describe along the way graphs whose MBD numbers are equal to $2$. This fills a gap in the literature. In fact, graphs $G$ with $\gmb(G) = 2$ were already described as follows. For this sake note that if $\gamma(G) = 1$, then $\gmb(G) = 1$, and if $\gamma(G)\ge 3$, then $\gmb(G)\ge 3$. 

\begin{proposition} {\rm\cite[Proposition 3.2]{gledel-2020}}
\label{prop:MB(G)=2}
Let $G$ be a graph with $\gamma(G)=2$. Then $\gmb(G)=2$ if and only if $G$ has a vertex that lies in two $\gamma$-sets of $G$.
\end{proposition}

If $x\in V(G)$ and $S\subseteq V(G)$, then we say that $S$ is {\em $x$-free} if $x\notin S$. We next complement Proposition~\ref{prop:MB(G)=2} with the following result. 

\begin{proposition} 
\label{prop:MB'(G)=2}
Let $G$ be a graph with $\gamma(G) = 2$. Then $\gmb'(G) = 2$ if and only if for any vertex $x\in V(G)$ there exists a vertex $y_x\in V(G)$ such that $y_x$ lies in two $x$-free $\gamma$-sets of $G$. 
\end{proposition}

\proof
Assume first that $\gmb'(G) = 2$. Let $x$ be an arbitrary vertex of $G$. Then $x$ is a legal first move of Staller. Since $\gmb'(G) = 2$, Dominator has an optimal reply $y$ such that he needs only one more move to finish the game. If $y'$ is such a vertex, then $\{y,y'\}$ is a $\gamma$-set of $G$. Now Staller can reply by selecting $y'$ as her second move. As $\gmb'(G) = 2$, Dominator can finish the game by selecting another vertex $y''$. Then $\{y,y'\}$ and $\{y,y''\}$ are two $\gamma$-sets of $G$ which are clearly $x$-free, hence setting $y_X = y$ we have the vertex required. 

Conversely, assume that for any vertex $x\in V(G)$ there exists a vertex $y_x\in V(G)$ such that $y_x$ lies in two $x$-free $\gamma$-sets of $G$. Assume that Staller starts the game by playing an arbitrary vertex $x$ of $G$. Then Dominator replies by selecting $y_x$ as her first move. (As $\gamma(G) = 2$, the vertex $y_x$ is clearly not a dominating vertex.)  Since $y_x$ lies in two $\gamma$-sets of $G$ (of cardinality $2$), Dominator can win the game in his second move. 
\qed

It is straightforward to check that in the case of corona products $G\odot H$, each of Propositions~\ref{prop:MB(G)=2} and~\ref{prop:MB'(G)=2} yields the same condition on $H$, that is, we have: 

\begin{corollary}
If $G$ is a graph with $n(G)\ge 2$, then $\gmb(G\odot H) = \gmb'(G\odot H)=2$ if and only if $n(G)=2$ and $H$ has a dominating vertex.
\end{corollary}

Graphs $G$ with $\gsmb(G)=2$ and with $\gsmb'(G)=2$ can be respectively described as follows. 

\begin{proposition} 
\label{prop:smb=2}
Let $G$ be a graph with $\delta(G)\geq 1$. Then the following hold.
\begin{enumerate}
\item[(i)] $\gsmb(G)=2$ if and only if $G$ has at least two strong support vertices.
\item[(ii)] $\gsmb'(G)=2$ if and only if $G$ has a strong support vertex.
\end{enumerate}
\end{proposition}

\proof
(i) Let $G$ be a graph with $\delta(G)\geq 1$ and assume that $\gsmb(G)=2$. Suppose that $G$ has at most one strong support vertex. Dominator first selects this strong support vertex. Now, Staller replies to this move with a vertex $v$. If at most one leaf is adjacent to $v$, then Dominator selects this leaf as his next turn and Staller cannot finish the game in her next move. Clearly, Staller can finish the game in her next move only when either $G$ has isolated vertices or at least two leaves adjacent to  $v$. Since $\delta(G)\geq 1$, $G$ has no isolated vertices. Thus $G$ has at least two strong support vertices.

Conversely, assume that $G$ has at least two strong support vertices. In the first move Staller can select a strong support vertex such that none of its leaves has been played by Dominator in his first move. Then Staller can isolate one leaf in her second move, hence $\gsmb(G) \leq 2$. Since $\delta(G) \geq 1$, we also have $\gsmb(G) \geq 2$. Therefore $\gsmb(G)= 2$.

(ii) Assume that $u$ us a strong support vertex of $G$. Since $\delta(G)\geq 1$, we have $\gsmb'(G) \geq 1 + \delta(G) \geq 2$. If Staller first selects $u$ in the S-game, then she can win the game in her second move, hence we also have $\gsmb'(G) \le 2$.
      
Conversely, assume that $\gsmb'(G)=2$. Since $1+\delta(G)\leq \gsmb'(G) = 2$, we see that $G$ has a leaf. Staller first selects a vertex $u$ which is adjacent to a leaf for otherwise $\gsmb'(G)\ge 3$.  If there is only one leaf adjacent to $u$, then Dominator selects that leaf as his next move. Since $\gsmb'(G)=2$,  Staller can finish the game in her next move. This is possible only when $G$ has an isolated vertex or $u$ has an unplayed leaf. Since $\delta(G)\geq 1$, $G$ has no isolated vertices. Therefore $u$ is a strong support vertex.    
\qed

Proposition~\ref{prop:smb=2} has been independently established by Bujt\'as and Dokyeesun and is contained in~\cite{dokyeesun-2024}. 

It is again straightforward to check that in the case of corona products $G\odot H$, each of the conditions (i) and (ii) of Propositions~\ref{prop:smb=2} and~\ref{prop:MB'(G)=2} yields the same condition on $H$ as follows. 

\begin{corollary}
\label{cor: smb=2}
If $G$ is a graph with $n(G)\geq 2$, then $\gsmb(G\odot H) = \gsmb'(G\odot H)=2$ if and only if $H$ has at least two isolated vertices.
\end{corollary}

We next describe corona products with both SMB numbers equal to three. 

\begin{proposition}
\label{prop: smb=3}
Let $G$ be a graph with $n(G)\geq 2$. Then $\gsmb(G\odot H) = \gsmb'(G\odot H) = 3$ if and only if either $H$ has at least two strong support vertices, or $H$ has an isolated vertex and at least one strong support vertex. 
\end{proposition}

\proof
Recall that Theorem~\ref{thm: formula 3} asserts $\gsmb(G\odot H) = \gsmb'(G\odot H) = 3$ if and only if $\gsmb(H)=2$. If $\delta(H) \ge 1$, then Proposition~\ref{prop:smb=2}(i) implies that this is if and only if $H$ has at least two strong support vertices. On the other hand, if $H$ has an isolated vertex, we infer that it must at least one strong support vertex. 
\qed

In general, we can characterize the graphs $G$ with $\gsmb'(G) =3$  as follows. 

\begin{proposition}\label{smb'=3}
    Let $G$ be a graph with $\delta(G) \geq 1$.   Then $\gsmb'(G) =3$ if and only if $G$ satisfies the following conditions.
    \begin{enumerate}
\item[(i)] $G$ has no strong support vertices.
\item[(ii)] There exist a vertex $u \in V(G)$ such that $G-u$ has at least two strong support vertices, or $G-u$ has an isolated vertex and a strong support vertex.
\end{enumerate}
\end{proposition}

\proof
Let $G$ be a graph with $\delta(G) \geq 1$.

Assume first that $\gsmb'(G) = 3$. Then $G$ has no strong support vertices, for otherwise we would have $\gsmb'(G) \le 2$. Let $u$ be an optimal first move of Staller in the S-game. 

\medskip\noindent{\bf Claim}: After Staller selects $u$, the remaining game is equivalent to the D-game played on $G-u$.

\medskip\noindent
To prove the claim we need to show that $u$ is not the vertex which is undominated (by Dominator) at the end of the game. Assume first that $ \deg(u)=1$. If $v$ is the neighbor of $u$, then Dominator must select $v$ in his first move, for otherwise Staller would win the game in her second move which is not possible by our assumption $\gsmb'(G) = 3$. Hence the claim. Assume second that $\deg(u)\ge 2$ and let $v$ be a vertex with smallest degree among the vertices from $N(u)$. Let Dominator select $v$ as his first move. Since $G$ has no strong support vertex, each vertex from $N(u)\setminus \{v\}$ has degree at least $2$, hence Staller cannot finish the game in her second move. This implies that the move $v$ is an optimal move of Dominator which in turn proves the claim. 

\smallskip
We have thus seen that after Staller selects $u$, the remaining game is a D-game on $G-u$. As $\gsmb'(G) = 3$, we have $\gsmb(G-u) = 2$. In view of Proposition ~\ref{prop:smb=2} this is possible only when either $G-u$ has at least two strong support vertices or $G-u$ has exactly one isolated vertex and a strong support vertex.

Conversely, assume that $G$ is a graph which satisfies $(i)$ and $(ii)$. Let $u$ be a vertex that satisfies $(ii)$ and let $u$ be the first move of Staller in the S-game. If $G-u$ has an isolated vertex then Dominator selects this vertex for otherwise Staller could claim it in the next move to win the game. After that, Staller selects a strong support vertex and wins the game in her third move. Similarly, if $G-u$ has no isolated vertex, then $G-u$ has two strong support vertices and again Staller can win in her third move. Thus we conclude that $\gsmb'(G) = 3$.
\qed

We conclude the section by a description of corona products with $\gmb$ equal to three. 

\begin{proposition}
Let $G$ be a graph with $n(G) \ge 2$ and let $H$ be a graph. Then $\gmb(G\odot H) = 3$ if and only if either $n(G)=2$ and $\gmb(H)=2$, or $n(G)=3$ and $H$ has a dominating vertex.
\end{proposition}

\proof
If $\gmb(G\odot H)=3$, then clearly $n(G) \le 3$. Assume first $n(G) = 2$. In their first two moves, Dominator and Staller will respectively select the two vertices of $G$, no matter whether the D-game or the S-game is played. After that, $\gmb(G\odot H)=3$ holds if and only if $\gmb(H) = 2$. Similarly, if $n(G) = 3$, then Dominator can win if and only if $H$ has a dominating vertex. 
\qed

\section{Concluding remarks}

In this paper, we have explored the Maker-Breaker domination game on corona products $G\odot H$ of graphs. It is worth noting the following equivalent approach to it. Let $\widehat{H}$ be the graph obtained from $H$ by adding a new vertex adjacent to all the vertices of $H$, and let $k\widehat{H}$ denote the disjoint union of $k$ copies of $\widehat{H}$. Then the Maker-Breaker domination game played on the corona product  $G\odot H$ is equivalent to the Maker-Breaker domination game played on $n(G)\widehat{H}$.

Clearly, $\gamma(\widehat{H}) = 1$. Having in mind the above equivalent formulation of the Maker-Breaker domination game on corona products, it would be of interest to extend the results of this paper to generalized coronas which are defined just as coronas, except that the graphs associated to vertices of $G$ can be different. Generalized coronas were by now investigated from different angles, cf.~\cite{dettleff-2016, michalski-2022, saravan-2021}. 

For the  the Maker-Breaker domination game on corona products, we have identified some general features, but several directions remain open for further exploration. Let us point to some of them. 

In Theorem~\ref{thm:bounds 1} we have bounded  the values of $\gmb(G\odot H)$ and of $\gmb'(G\odot H)$ for the case when 
$o(H)=\mathcal{D}$. It would be interesting to investigate possible values that $\gmb'(G\odot H) - \gmb(G\odot H)$ can attain.

In Theorem~\ref{33} we have determined $\gmb(G\odot T)$ and $\gmb'(G\odot T)$ for an arbitrary graph $G$ with $n(G)\ge 2$ and an arbitrary tree $T$ with a perfect matching. It would be interesting to derive such a result also for arbitrary trees, that is, also for trees which do not contain perfect matchings. 

With respect to small MBD numbers it in particular remains to consider arbitrary graphs $G$ with $\gsmb(G) = 3$, or with $\gmb(G) = 3$, or with $\gmb'(G) = 3$, as well as corona products $G\odot H$ with $\gmb'(G\odot H) = 3$. 

\section*{Acknowledgements}

Sandi Klav\v zar was supported by the Slovenian Research Agency ARIS (research core funding P1-0297 and projects N1-0218, N1-0285). 

\section*{Declaration of interests}
 
The authors declare that they have no conflict of interest. 

\section*{Data availability}
 
Our manuscript has no associated data.

\end{document}